\documentclass[a4paper]{amsart}
\usepackage{enumitem, amsmath, amssymb, amsthm, xspace}
\usepackage[bookmarksnumbered,bookmarksopen]{hyperref}

\providecommand{\hyperref}[2][]{#2}
\providecommand{\hypertarget}[2]{#2}
\providecommand{\hyperlink}[2]{#2}

\title[Matrix types of ultramatricial algebras]{Characterization of matrix
  types of ultramatricial algebras}
\author{G\'abor Braun}
\address{Alfr\'ed R\'enyi Institute of Mathematics\\
  Hungarian Academy of Sciences\\
  Budapest\\
  Re\'altanoda u 13--15\\
  1053\\
  Hungary}
\date{}
\email{braung@renyi.hu}

\AtEndDocument{\small Published version available at \url{http://nyjm.albany.edu:8000/j/2005/11-2.html}.}

\newtheorem{theorem}{Theorem}[section]
\newtheorem{lemma}[theorem]{Lemma}
\newtheorem{proposition}[theorem]{Proposition}
\newtheorem{corollary}[theorem]{Corollary}

\theoremstyle{definition}
\newtheorem{definition}[theorem]{Definition}

\theoremstyle{remark}
\newtheorem{remark}[theorem]{Remark}

\newcommand{\setZ}{\ensuremath{\mathbb{Z}}}
\newcommand{\setQ}{\ensuremath{\mathbb{Q}}}

\newcommand*{\completion}[1]{\widehat{#1}}

\DeclareMathOperator{\Aut}{Aut}
\DeclareMathOperator{\mt}{mt}

\newcommand*{\nSeq}[1][n]{\ensuremath{{\csname HCode\endcsname{}}\sp{#1}2}\xspace}\newcommand{\InfSeq}{\ensuremath{{\csname HCode\endcsname{}}\sp{<\omega}2}\xspace}\newcommand{\AtMostnSeq}{\ensuremath{{\csname HCode\endcsname{}}\sp{\leq n}2}\xspace}
\begin{document}
\begin{abstract}
  For any equivalence relation \(\equiv\) on positive integers such that \(nk \equiv mk\) if and
  only if \(n \equiv m\), there is an abelian group \(G\) such that the
  endomorphism ring of \(G^n\) and \(G^m\) are isomorphic if and only if \(n \equiv
  m\).  However, \(G^n\) and \(G^m\) are not isomorphic if \(n \neq m\).
\end{abstract}

\maketitle{}

\section{Introduction}
\label{sec:introduction}

We construct partially ordered abelian groups such that the orbit of a
distinguished element under the automorphism group is prescribed; the precise
statement is Theorem~\ref{th:2} in Subsection~\ref{sec:dimension-groups}.  The
prescribed orbit controls the \emph{matrix type} of a ring, i.e., which matrix
algebras over the ring are isomorphic, hence we can characterize the matrix
types of ultramatricial algebras over any principal ideal domain, see
Theorem~\ref{th:1} in Subsection~\ref{sec:ultr-algebr}.  If the ground ring is
\(\setZ\) then these algebras are realizable as endomorphism rings of torsion-free
abelian groups, which groups therefore have the property stated in the
abstract, see Corollary~\ref{cor:1}.

We are indebted to P\'eter V\'amos who draw our attention to this wonderful problem.

\subsection{Dimension groups}
\label{sec:dimension-groups}

An \emph{order unit} in a partially ordered abelian group is a positive
element \(u\) such that for every element \(x\) there is a positive integer
\(n\) such that \(nu \geq x\). A \emph{dimension group} \((D, \leq, u)\) is a
countable partially ordered abelian group \((D, \leq)\) with order unit \(u\)
such that every finite subset of \(D\) is contained in a subgroup, which is
isomorphic to a direct product of finitely many copies of \((\setZ, \leq)\) as a
partially ordered abelian group.

Our main result, which will be proven from Section~\ref{sec:overv-constr} on, is:
\begin{theorem}
  \label{th:2}
  Let \(H \leq \setQ_+^\times\) be a subgroup of the multiplicative group of the positive
  rational numbers. Then there is a dimension group \((D, \leq, u)\) whose group
  of order-preserving automorphisms is isomorphic to \(H\). Furthermore, under
  this isomorphism every element of \(H\) acts on \(u\) by multiplication by
  itself as a rational number.
\end{theorem}

In the special case when \(H\) is generated by a set \(S\) of prime numbers,
one may choose \(D\) to be the ring \(\setZ[S^{-1}]\) and \(u=1\), see
\cite[Proposition~4.2]{MR99f:16008}.

\subsection{Ultramatricial algebras}
\label{sec:ultr-algebr}

An \emph{ultramatricial algebra} over a field or principal ideal domain \(F\)
is an \(F\)-algebra which is a union of an upward directed countable set of
\(F\)-subalgebras, which are direct products of finitely many matrix algebras
over \(F\).

\subsubsection{Matrix types of rings}
\label{sec:matrix-types-rings}

Let \(M_n(R)\) denote the ring of \(n \times n\) matrices over the ring \(R\).
Obviously, \(M_n(R)\) is ultramatricial if \(R\) is ultramatricial.  The
\emph{matrix type} of a ring \(R\) is the equivalence on positive integers
describing which matrix algebras over \(R\) are isomorphic:
\begin{equation}
  \label{eq:1}
  \mt(R) := \{ (n,m) \mid M_n(R) \cong M_m(R) \}.
\end{equation}
Clearly, if \((n,m) \in \mt(R)\) then \((mk,nk) \in \mt(R)\) for all positive
integers \(m\), \(n\), \(k\).  The converse also holds for ultramatricial
algebras but probably not for all rings.

The next theorem states that all such equivalences indeed arise as matrix
types of ultramatricial algebras:
\begin{theorem}
  \label{th:1}
  Let \(F\) be a field or principal ideal domain and \(\equiv\) be an equivalence
  relation on the set of positive integers. Then the following are equivalent:
  \begin{enumerate}[label=(\roman*)]
  \item For all positive integers \(n\), \(m\) and \(k\)
    \begin{equation}
      \label{eq:2}
      n \equiv m \iff nk \equiv mk.
    \end{equation}
  \item There is a (unique) subgroup \(H\) of the multiplicative group
    \(\setQ_+^\times\) of positive rational numbers such that for all positive integers
    \(n\) and \(m\):
    \begin{equation}
      \label{eq:29}
      n \equiv m \iff \frac{n}{m} \in H.
    \end{equation}
  \item There exists an ultramatricial \(F\)-algebra with matrix type \(\equiv\).
  \end{enumerate}
\end{theorem}

The second statement is clearly a reformulation of the first one, which is
useful for explicit construction of equivalences, as pointed out by the
referee.

The equivalence of the second and third statements is a simple consequence of
our main result, Theorem~\ref{th:2}, as we will explain in the
\hyperref[sec:equiv-dimens-groups]{next section}.

\subsubsection{Basis types of rings}
\label{sec:basis-types-rings}

In a similar vein, the \emph{basis type} of a ring \(R\) characterizes which
finite rank free modules are isomorphic:
\begin{equation}
  \label{eq:30}
  \operatorname{bt}(R) := \{ (n,m) | R^n \cong R^m \}.
\end{equation}
The analogue of Theorem~\ref{th:1} for basis types is Theorem~1 of
\cite{MR0238903}: every equivalence relation with the property \(n \equiv m\) if
and only if \(n+k \equiv m+k\) is the basis type of a ring of infinite matrices.

Clearly, if \(R^n \cong R^m\) then \(R^{n+k} \cong R^{m+k}\).  All equivalence
relation with this property arise as basis type of a ring: see
\cite{MR0132764}.

The basis type is obviously smaller than the matrix type.  It seems plausible
that this is the only relation between the two types.  For example,
ultramatricial algebras have trivial basis type i.e., free modules of
different finite rank are not isomorphic.

Every ultramatricial algebra \(R\) over \(F=\setZ\) is a countable reduced
torsion-free ring, and hence is the endomorphism ring of a torsion-free
abelian group \(G\) by \cite[Theorem~A]{MR27:3704}.  Then \(M_n(R)\) is the
endomorphism ring of \(G^n\).  Obviously \(G^n\) and \(G^m\) are not
isomorphic if \(n \neq m\) since there is no invertible \(n \times m\) matrix over
\(R\).  Hence we have the statement in the abstract as an immediate corollary
to the last theorem:
\begin{corollary}\label{cor:1}
  Let \(\equiv\) be an equivalence relation on positive integers with the property
  that \(n \equiv m\) if and only if \(nk \equiv mk\) for all positive integers \(n\), \(m\)
  and \(k\).  Then there is a torsion-free abelian group \(G\) such that the
  endomorphism ring of \(G^n\) and \(G^m\) are isomorphic if and only if
  \(n \equiv m\).  Moreover, \(G^n\) and \(G^m\) are isomorphic if and only if
  \(n=m\).
\end{corollary}

\section{Equivalence of dimension groups and ultramatricial algebras}
\label{sec:equiv-dimens-groups}

Dimension groups and ultramatricial algebras over a fixed field or principal
ideal domain are essentially the same.  In this section, we recall this
equivalence, which shows that Theorem~\ref{th:1} follows from
Theorem~\ref{th:2}.  For more details and proofs see
\cite[Proposition~4.1]{MR99f:16008} or \cite[Chapter~15, Lemma~15.23,
Theorems~15.24 and 15.25]{MR93m:16006}, which assume that the ground ring is a
field, but the arguments also work when it is a principal ideal domain.

First we define the functor \(K_0\) from the category of rings to the category
of preordered abelian groups with a distinguished order unit.  The isomorphism
classes of finitely generated projective left modules over a ring \(R\) form a
commutative monoid where the binary operation is the direct sum.  The quotient
group of the monoid is denoted by \(K_0(R)\).  Declaring the isomorphism
classes of projective modules to be nonnegative, \(K_0(R)\) becomes a
preordered group.  The isomorphism class \(u\) of \(R\), the free module of
rank \(1\), is an order unit of \(K_0(R)\).

If \(f\colon R \to S\) is a ring homomorphism then let \(K_0(f)\) map the isomorphism
class of a projective module \(P\) to that of \(S \otimes_R P\).  So \(K_0(f)\) is a
homomorphism preserving both the order and the order unit.

If \(R\) is an ultramatricial algebra then \(K_0(R)\) is a dimension group.
Conversely, every dimension group is isomorphic to the \(K_0\) of an
ultramatricial algebra.  If \(R\) and \(S\) are ultramatricial algebras then
every morphism between \(K_0(R)\) and \(K_0(S)\) is of the form \(K_0(f)\) for
some algebra homomorphism \(f\colon R \to S\).  However, \(f\) is not unique in
general.  Nevertheless, every isomorphism between \(K_0(R)\) and \(K_0(S)\)
comes from an isomorphism between \(R\) and \(S\).

Thus, by restriction, \(K_0\) is essentially an equivalence between the category of
ultramatricial algebras over a given field and the category of dimension
groups with morphisms the group homomorphisms preserving both the order and
the order unit. 

Now we examine how the matrix type of an ultramatricial algebra can be
recovered from its dimension group.  The standard Morita equivalence between
\(R\) and \(M_n(R)\) induces an isomorphism between the \(K_0\) groups.
However, this isomorphism does not preserve the order unit, in fact
\(K_0(M_n(R))\) is \((K_0(R), \leq, nu)\) where \(u\) is the order unit of
\(K_0(R)\).

So if \(R\) is an ultramatricial algebra, then the dimension groups
\(K_0(M_n(R))\) and \(K_0(M_m(R))\) (and hence the algebras \(M_n(R)\) and
\(M_m(R)\)) are isomorphic if and only if there is an order-preserving
automorphism of \(K_0(R)\) sending \(nu\) to \(mu\).  Obviously, for any
dimension group \((D, \leq, u)\) whether an order-preserving automorphism maps
\(nu\) to \(mu\), depends only on the factor \(m/n\).  Such factors \(m/n\)
form a subgroup of the multiplicative group \(\setQ_+^\times\) of the positive
rationals.

So the classification of matrix types of ultramatricial algebras is equivalent
to the classification of subgroups of \(\setQ_+^\times\) which arise from dimension
groups in the above construction. Theorem~\ref{th:2} states that all subgroups
arise, and Theorem~\ref{th:1} is just the translation of it to the language of
matrix types of ultramatricial algebras.

\section{Overview of the construction}
\label{sec:overv-constr}

In the rest of the paper we prove Theorem~\ref{th:2}.  In this section we
outline the main ideas of the proof and leave the details for the following
sections.  The \hyperref[sec:notation]{next section} fixes notations used frequently in the rest of the
paper.  Section~\ref{sec:abelian-groups-endo} recalls a construction of
abelian groups.  The actual proof is contained in the rest of the sections,
which are organized so that they can be read independently.  At the beginning
of every section, we shall refer to its main proposition, which will be the
only statement used in other sections.  The same is true for subsections.

To prove Theorem~\ref{th:2}, we fix a subgroup \(H\) of the positive rationals
and construct a dimension group \(D\) for it.  We search for \(D\) (as an abelian
group without any order) in the form \(D := \setQ u \oplus G\) where \(H\) acts on the
direct sum componentwise.  We let \(H\) act on \(\setQ u\) by multiplication as
required to act on the order unit.  A key observation
(Proposition~\ref{prop:2}) is that if the only automorphisms of \(G\) are the
elements of \(H\) and their negatives, then for any partial order \(\leq\) on
\(D\), which is preserved by \(H\) and makes \((D, \leq, u)\) a dimension group,
the order-preserving automorphism group of \(D\) is only \(H\).  So \((D, \leq,
u)\) satisfies the \hyperref[th:2]{theorem}.

Therefore, all we have to do is to find such a \(G\) and \(\leq\).  Actually,
\(G\) is already constructed by A.~L.~S.~Corner in~\cite{MR27:3704}.  Since we
shall use the structure of \(G\) to define the partial order \(\leq\), we recall
a special case of the construction in Section~\ref{sec:abelian-groups-endo},
which is enough for our purposes.  See also \cite{GT} for the general
statement.

Finally, we define a partial order \(\leq\) on \(D\) making it a dimension group
in Subsection~\ref{sec:partial-order}.  The basic idea is to explicitly make
some subgroups of \(D\) order-isomorphic to \(\setZ^n\) and show that these
partial orders on the subgroups are compatible.  This will be based on a
description of \(G\) (Proposition~\ref{prop:3}): we define elements of \(G\),
which will be bases of order-subgroups \(\setZ^n\) of \(D\) and state relations
between them implying the compatibility of partial orders.

\section{Notation}
\label{sec:notation}

This section is just for fixing notation used throughout the paper.

Let \(x_h\) denote the element of the group ring \(\setZ H\) corresponding to \(h
\in H\). This is to distinguish \(x_ha\) the value of \(x_h\) at \(a\) from
\(ha\) the element \(a\) multiplied by the rational number \(h\).

If the automorphism group of an abelian group is the direct product of a group
\(H\) and the two element group generated by \(-1\) then we say that the
automorphism group is \(\pm H\).

\section{Abelian groups with prescribed endomorphisms}
\label{sec:abelian-groups-endo}

In this section we revise a special case of A.~L.~S.~Corner's construction of
abelian groups with prescribed endomorphism rings.  I am grateful to R\"udiger
G\"obel and his group for teaching me this method.

Let \(\completion{M}\) denote the \(\setZ\)-adic completion of the abelian group \(M\).
The following result is is a special case of Theorem~1.1 from
\cite{0178.02303} by taking \(A=R\) and \(N_k=0\):

\begin{proposition}
  \label{prop:1}
  Let \(R\) be a ring with free additive group. Let \(B\) be a free
  \(R\)-module of rank at least \(2\) and \(\{w_b : b \in B \setminus \{0\}\}\) a collection
  of elements of \(\completion{\setZ}\) algebraically independent over \(\setZ\). Let the
  \(R\)-module \(G\) be
  \begin{equation}
    \label{eq:3}
    G := \langle B, Rbw_b : b \in B \setminus \{0\} \rangle_* \subseteq \completion{B}.
  \end{equation}
  
  Then \(G\) is a reduced abelian group with endomorphism ring \(R\).
\end{proposition}
Here \(*\) means purification: i.e., we add all the elements \(x\) of
\(\completion{B}\) to \(E := B \oplus_{b \in B \setminus \{0\}} Rbw_b\) for which \(nx \in E\)
for some nonzero integer \(n\).  The usual down-to-earth description of \(G\)
is the following, which we shall use in
Subsection~\ref{sec:abelian-group-under}: we select positive integers \(m_n\)
such that every integer divides \(m_1 \dotsm m_n\) for \(n\) large enough.  We
choose elements \(w_b^{(n)}\) of \(\completion{\setZ}\) for all nonzero elements
\(b\) of \(B\) and natural numbers \(n\) such that \(w_b^{(0)} = w_b\) and
\(w_b^{(n)} - m_n w_b^{(n+1)}\) is an integer.  Then \(G\) is generated by the
submodules \(B\) and \(Rbw_b^{(n)}\).  Obviously, \(G_n := B \oplus_{b \in B \setminus \{0\}}
Rbw_b^{(n)}\) is a submodule of \(G\) and these modules \(G_n\) form an
increasing chain whose union is the whole \(G\).  It is important to note that
\(G\) does \emph{not} depend on the choice of the \(m_n\) and \(w_b^{(n)}\).
However, the submodules \(G_n\) \emph{does} depend on them.

\begin{remark}
  \label{remark:Corner}
  Note that the automorphism group of \(G\) is the group of units of \(R\),
  which is just \(\pm H\) if \(R = \setZ H\) and \(H\) is an orderable group.  (It
  is a famous conjecture that the group of units of \(\setZ H\) is \(\pm H\) for all
  torsion-free groups \(H\).)
  
  We will be interested in the case when \(H \leq \setQ_+^\times\) and \(R = \setZ H\) is a
  group ring.  The free module \(B\) will have countable rank.  Recall that
  there are continuum many elements of \(\completion{\setZ}\) algebraically independent
  over \(\setZ\), so the construction works in this case, and we will get a reduced
  abelian group \(G\) with automorphism group \(\pm H\).
\end{remark}

\section{Construction of the dimension group}
\label{sec:constr-dimens-group}

We now construct our dimension group \(D = \setQ u \oplus G\) starting from a subgroup
\(H \leq \setQ_+^\times\) of the positive rationals.  The main result of the section is
Proposition~\ref{prop:4}.

The \hyperref[prop:3]{proposition} in Subsection~\ref{sec:abelian-group-under}
defines \(G\) and summarizes the (technical) properties of \(G\) used in
Subsection~\ref{sec:partial-order} to put the partial order on \(D\).

Let \(\nSeq\) denote the set of sequences of length \(n\) whose elements are
\(0\) or \(1\). If \(y\) is a finite sequence of \(0\) and \(1\) then let
\(y0\) denote the sequence obtained from \(y\) by adding an additional element
\(0\) to the end. Similarly, we can define \(y1\).  These sequences will be
identified with elements of \(G\).

\subsection{The abelian group under the dimension group}
\label{sec:abelian-group-under}

In this subsection we construct an abelian group \(G\) satisfying the
following properties. Essentially, \(G\) will be the underlying abelian group
of our dimension group.

\begin{proposition}\label{prop:3}
  Let \(H \leq \setQ_+^\times\) be a subgroup of the multiplicative group of the positive
  rational numbers. Let \(R:=\setZ H\) denote its group ring. Then there exist:
  \begin{itemize}
  \item an \(R\)-module \(G\), which is also a reduced abelian group,
  \item a finite subset \(F_n \subseteq H\) for all positive integer \(n\),
  \item positive integers \(s_n\), \(t_n\), \(k_n\) and \(l_n\) for \(n \geq 0\)
  \end{itemize}
  subject to the following conditions:
  \begin{enumerate}[label=(\roman*)]
  \item \label{item:1} \(\Aut G = \pm H\).
  \item \label{item:11} Structure of \(G\):
    \begin{enumerate}[label=(\alph*), ref=\theenumi.(\alph*)]
    \item \label{item:2} \(G\) is a union of an increasing sequence of free
      submodules \(G_n\).
    \item \label{item:3} \(G_n\) has base \(\nSeq \cup \{ c_i^{(n)} : i = 1, \dotsc,
      n \}\).
    \end{enumerate}
  \item \label{item:12} Relations describing the inclusion \(G_n \subseteq G_{n+1}\):
    \begin{enumerate}[label=(\alph*), ref=\theenumi.(\alph*)]
    \item \label{item:4} \(y = y0 + s_n^2 \cdot y1\) \quad for all \(y \in \nSeq\).
    \item \label{item:5} There are integers \(n_{h,y}^{(i)}\) for all \(i \leq
      n\), \(h \in H\) and \(y \in \nSeq\) such that
      \begin{equation}\label{eq:12}
        c_i^{(n)} - s_nt_n c_i^{(n+1)} = \sum_{\substack{h \in F_i \\ y \in \nSeq}}
        n_{h,y}^{(i)} x_h \cdot y0, \quad 0 \leq n_{h,y}^{(i)} < s_nt_n.
      \end{equation}
    \end{enumerate}
  \item \label{item:10} Properties of \(s_n\), \(t_n\), \(k_n\) and \(l_n\):
    \begin{enumerate}[label=(\alph*), ref=\theenumi.(\alph*)]
    \item \label{item:7}
      \begin{align*}
        k_0 &= 1, & k_{n+1} &= s_n k_n,\\
        l_0 &= 1, & l_{n+1} &= t_n l_n.
      \end{align*}
    \item \label{item:8} Every positive integer divides \(k_n\) and \(l_n\) for
      \(n\) large enough.
    \item \label{item:9} \(k_n (s_n^2-1) \geq l_n s_n t_n \sum_{i=1}^n \sum_{h \in F_i} h\).
    \end{enumerate}
  \end{enumerate}
\end{proposition}

\begin{proof}
The construction of the items is easy.  One has to care about defining
them in the correct order.

Let \(B\) be a countable-rank free module over the group ring \(R:=\setZ H\) i.e.,
\(B = \setZ H \otimes A\) where \(A\) is a free abelian group of countably infinite
rank.  Using this group we define \(G\) by Equation~\eqref{eq:3} as in
Proposition~\ref{prop:1}.  (The \(\setZ\)-adic integers \(w_b\) can be chosen
arbitrarily.)  The \hyperref[prop:1]{proposition} tells us that \(G\) is a
reduced abelian group and \(\Aut G = \pm H\) so \ref{item:1} is satisfied.

We identify a base of \(A\) with the finite sequences of \(0\) and \(1\) not ending with \(0\)
(so with the sequences ending with \(1\) and the empty sequence). In this
base, the sequences of length at most \(n\) is a basis of a free subgroup
\(A_n\) of \(A\) and the free \(R\)-module \(B_n = R \otimes A_n\).

Let us enumerate the elements of \(B \setminus \{0\}\) into a sequence \(b_1, b_2,
\dotsc\) such that \(b_n \in B_n\). Every element of \(B\) can be written uniquely
as \(\sum_{h \in H} x_h b_h\) where \(b_h \in A\) and only finitely many of the
\(b_h\) are nonzero. We define the \emph{support} of an element of \(B\) as the
finite set
\begin{equation}
  \label{eq:14}
  \left[ \sum_{h \in H} x_h b_h \right] := \{ h \in H \mid b_h \neq 0 \} \quad (b_h \in A).
\end{equation}
Let \(F_i := [b_i]\) be the support of \(b_i\).

Now we are ready to define our positive integers \(s_n\), \(t_n\), \(k_n\) and
\(l_n\). Let us impose the following additional condition on them:
\begin{trivlist}
\item[\hypertarget{item:divisible}{(*)}] \(t_n\) is divisible by \(n\), and \(t_n\) divides \(s_n\).
\end{trivlist}
Now the integers can be defined recursively such that \ref{item:7}, \ref{item:9}
and~\hyperlink{item:divisible}{(*)} hold. These automatically imply the truth of \ref{item:8}.

Now we identify the sequences of \(0\) and \(1\) with elements of \(A\).  We
have already done this for the sequences not ending with zero: they form a
basis of \(A\).  As dictated by~\ref{item:4}, we set
\begin{equation}
  \label{eq:15}
  y0 := y - s_n^2 \cdot y1 \quad (y \in \nSeq).
\end{equation}
(This is in fact a recursive definition on the length of \(y\) since \(y\) may
also end with zero.) Thus the sequences of finite length are identified with
elements of \(A\) such that \ref{item:4} holds and the sequences of length
exactly \(n\) form a basis of \(A_n\) as can be easily seen by induction on
\(n\).

We turn to the definition of the \(G_n\). Let us choose \(\setZ\)-adic integers
\(w_i^{(n)}\) for \(1 \leq i \leq n\) such that
\begin{equation}
  \label{eq:16}
  w_i^{(0)} := w_{b_i}, \quad w_i^{(n)} - s_n t_n w_i^{(n+1)} \in \setZ, \quad
  w_i^{(n)} \in \completion{\setZ}.
\end{equation}
We let \(G_n\) be the free submodule
\begin{equation}
  \label{eq:17}
  G_n := B_n \oplus \bigoplus_{i=1}^n R b_i w_i^{(n)}.
\end{equation}
It follows from the definition of \(G\) (Equation~\eqref{eq:3}) that the groups
\(G_n\) form an increasing sequence of submodules whose union is \(G\).

The only missing entities are the elements \(c_i^{(n)}\).  We could set \(c_i^{(n)}=b_i
w_i^{(n)}\) to satisfy \ref{item:3} but this may not be appropriate
for~\ref{item:5}.  Therefore we shall set
\begin{equation}
  \label{eq:18}
  c_i^{(n)} := b_i w_i^{(n)} + b_i^{(n)} \quad (b_i^{(n)} \in B_n, i \leq n)
\end{equation}
for some \(b_i^{(n)}\).  This ensures \ref{item:3}.  For \(i\) fixed, we are
going to define the \(b_i^{(n)}\) recursively for \(n \geq i\) subject to:
\begin{enumerate}[label=(\Alph*)]
\item \label{item:13} \(b_i^{(n)} \in B_n\).
\item \label{item:14} \([b_i^{(n)}] \subseteq [b_i] = F_i\).
\item \label{item:15} For suitable integers \(n_{h,y}^{(i)}\):
  \begin{equation}
    \label{eq:19}
    b_i \left(w_i^{(n)} - s_n t_n w_i^{(n+1)} \right) + b_i^{(n)} =
    s_n t_n b_i^{(n+1)} + \sum_{h \in F, y \in \nSeq} n_{h,y}^{(i)} x_h \cdot y0, \quad 0 \leq n_{h,y}^{(i)} < s_nt_n.
  \end{equation}
\end{enumerate}
Note that the last equation is just a reformulation of~\ref{item:5} in terms of
the \(b_i^{(n)}\).

Now we carry out the recursive definition. We can start with
\(b_i^{(i)}:=0\). Observe that \ref{item:15} determines how to define
\(b_i^{(n+1)}\): the left-hand side is an element of \(B_{n+1}\), a free
abelian group with basis \(x_h y\) for \(h \in H\) and \(y\) a sequence of \(0\)
and \(1\) of length \(n+1\). We divide the coefficient of every \(x_h y\) by
\(s_n t_n\). The quotient gives the coefficient of \(x_h y\) in
\(b_i^{(n+1)}\) and the remainder is a coefficient of the big sum on the
right. Since the support of the left-hand side is contained in \(F_i\) by
induction, the same is true for \(b_i^{(n+1)}\) and the sum on the right-hand
side.  The left-hand side is actually contained in \(B_n\) not only
\(B_{n+1}\). This means that the coefficients of sequences of length \(n+1\)
ending with \(1\) are divisible by \(s_n^2\) by~\ref{item:4} and hence by
\(s_n t_n\) since \(t_n\) divides \(s_n\). So in the sum on the right-hand
side, the coefficient of sequences ending with \(1\) is zero. Thus we have
defined \(b_i^{(n+1)}\) according to the requirements. 
\end{proof}

\subsection{The partial order}
\label{sec:partial-order}

In this subsection we define a partial order on \(D\) which will make it a
dimension group.

\begin{proposition}\label{prop:4}
  Let \(H \leq \setQ_+^\times\) be a subgroup of the multiplicative group of the positive
  rational numbers acting on \(\setQ u\) by multiplication. Suppose \(G\) is a
  group satisfying the conditions of Proposition~\ref{prop:3}. Let \(H\) act
  on \(D := \setQ u \oplus G\) componentwise. Then there is a partial order \(\leq\) on
  \(D\) such that \((D, \leq, u)\) is a dimension group on which \(H\) acts by
  order-preserving automorphisms.
\end{proposition}

The dimension group \(D\) in the \hyperref[prop:4]{proposition} satisfies all
requirements of Theorem~\ref{th:2}.  We shall see in the
\hyperref[sec:autom-dimens-groups]{next section} that the group of
order-preserving automorphisms of \(D\) is exactly \(H\). The other
requirements are obviously satisfied.

\begin{proof}
We define the partial order on a larger group, the divisible hull \(\setQ D\) of
\(D\).  For all natural number \(n\) and finite subset \(F\) of \(H\) we
define a subgroup of \(D\):
\begin{equation}
  \label{eq:5}
  D_{n,F} := \setZ \frac{u}{k_n} \oplus \bigoplus_{\substack{h \in F \\ y \in \nSeq}} \setZ x_hy \oplus
  \bigoplus_{\substack{i=1 \\ h \in F}}^n \setZ x_hc_i^{(n)}.
\end{equation}
We define the partial order on \(\setQ D_{n,F}\) as the product order
\begin{equation}
  \label{eq:6}
  (\setQ D_{n,F}, \leq) := (\setQ v_{n,F}, \leq) \times \prod_{\substack{h \in F \\ y \in \nSeq}} (\setQ x_hy, \leq) \times
  \prod_{\substack{i=1 \\ h \in F}}^n (\setQ x_hc_i^{(n)}, \leq)
\end{equation}
where
\begin{equation}
  \label{eq:7}
  v_{n,F} := \frac{u}{k_n} - \sum_{h \in F} h^{-1} x_h
  \left(
    k_n \sum_{y \in \nSeq} y +
    l_n \sum_{i=1}^n c_i^{(n)}
  \right).
\end{equation}
Note that this makes \(u\) an order unit of \(\setQ D_{n,F}\).

The subgroups \(D_{n,F}\) form a directed system whose union is the whole
\(D\). It follows that the subgroups \(\setQ D_{n,F}\) form a directed system
whose union is \(\setQ D\).

We are in a position now to reduce the proof to two lemmas stated below.  The
\hyperref[lem:3]{first one} states that the inclusions between the \(\setQ
D_{n,F}\) are order-embeddings.  The \hyperref[lem:4]{second one} claims that
a cofinal set of the \(D_{n,F}\) is order-isomorphic to a direct product of
finitely many copies of \((\setZ, \leq)\).

It follows that \(\setQ D\) has a unique partial order which extends the partial
order of all the \(\setQ D_{n,F}\).  Under this partial order \(u\) is clearly an
order-unit.  The \(D_{n,F}\) provide enough order-subgroups isomorphic to a
finite direct power of \((\setZ, \leq)\), hence \((D, \leq, u)\) is a dimension group.

The partial order is preserved by \(H\) since any \(h \in H\) maps \(\setQ D_{n,F}\)
bijectively onto \(\setQ D_{n,hF}\) and this bijection is an order-isomorphism of
the two subgroups. (Note that \(x_hv_{n,F} = h v_{n,hF}\).)

Thus \((D, \leq, u)\) have all the properties claimed.

All in all, the proposition is proved modulo the following two lemmas.
\end{proof}

\begin{lemma}
  \label{lem:3}
  The inclusions between the subgroups \(\setQ D_{n,F}\) are order-embeddings.
\end{lemma}

\begin{lemma}
  \label{lem:4}
  A cofinal set of the groups \(D_{n,F}\) is order-isomorphic to a finite direct
  power of \((\setZ, \leq)\).
\end{lemma}

We consider first the inclusions.

\begin{proof}[Proof of Lemma~\ref{lem:3}]
First we prove the claim that \(\setQ D_{n,F}\) is an order-subgroup of \(\setQ
D_{n+1,F'}\) if \(F'\) contains \(F\) and \(FF_i\) for all \(i \leq n\).
For this, it is good to have the following general example of an
order-embedding of \((\setQ, \leq)^m\) into \((\setQ, \leq)^{m+k}\) given by a matrix:
\begin{equation}
  \label{eq:28}
  \begin{pmatrix}
        >0 & 0  & \ldots & 0               & \geq0 & \ldots & \geq0\\
    0  & >0 & 0  & \vdots               & \geq0 & \ldots & \geq0\\
    \vdots  &    & \ddots  & \vdots               & \vdots  &  \vdots  & \vdots \\
    0  & \ldots & 0  & >0              & \geq0 & \ldots & \geq0\\
  \end{pmatrix}.
\end{equation}
On the first \(m\) coordinates this is an order-isomorphism: every coordinate
is multiplied by a positive number. On the last \(k\) coordinates the map is
an arbitrary order-preserving map. By permutating the coordinates, we may
complicate the map.

All in all, a homomorphism \((\setQ, \leq)^m \to (\setQ, \leq)^M\) is an order-embedding if
the canonical basis elements of the domain have only nonnegative coordinates
in the codomain and every basis element has a positive coordinate, which
coordinate is zero for the other basis elements.

We show that the inclusion of \(\setQ D_{n,F}\) into \(\setQ D_{n+1,F'}\) is an
order-embedding of the above type using the direct product
decomposition~\eqref{eq:6}.  To this end, we express the generators of \(\setQ
D_{n,F}\) as linear combination of the generators of \(\setQ D_{n+1,F'}\) (which
in particular shows that \(\setQ D_{n,F}\) is really a subgroup of \(\setQ
D_{n+1,F'}\)):
\begin{align}
  \label{eq:8}
  x_h y &= x_h \cdot y0 + s_n^2 x_h \cdot y1 & y \in \nSeq,\, &h \in F\\
  \label{eq:9}
  x_h c_i^{(n)} &= s_n t_n x_h c_i^{(n+1)} + \sum_{\substack{t \in F_i \\ y \in
      \nSeq}} n_{t,y}^{(i)} x_{ht} \cdot y0 & i \leq n,\, &h \in F\\
  \label{eq:10}
  v_{n,F} &=
  \begin{aligned}[t]
    s_n v_{n+1,F'} &+ \sum_{h \in F'} s_n l_{n+1} h^{-1} x_h c_{n+1}^{(n+1)}\\
    &+ \sum_{\substack{h \in F' \setminus F \\ i \leq n}} s_n l_{n+1} h^{-1} x_h c_i^{(n+1)}\\
    &+ \sum_{\substack{h \in F \\ y \in \nSeq}} k_n (s_n^2 - 1) h^{-1} x_h \cdot y0\\
    &+ \sum_{\substack{h \in F' \setminus F \\ y \in \nSeq[n+1]}} k_n s_n^2 h^{-1} x_h y\\
    &- \sum_{\substack{h \in F, i \leq n \\t \in F_i, y \in \nSeq}} l_n n_{t,y}^{(i)}
    h^{-1} x_{ht} \cdot y0
  \end{aligned}
\end{align}
These are easy consequences of the formulas under~\ref{item:12} of
Proposition~\ref{prop:3} and the definition~\eqref{eq:7} of \(v_{n,F}\).  All
the coordinates of the above generators are obviously nonnegative except for
the coefficient of \(x_h y0\) of \(v_{n,F}\) for \(y \in \nSeq\) and \(h \in F'\).
So let us consider the coefficient of \(h^{-1} x_h y0\) in
Equation~\eqref{eq:10}: from the third or forth row comes \(k_n (s_n^2-1)\)
or \(k_n s_n^2\) depending on whether \(h\) is contained in \(F\).  From the
last row \(- l_n n_{t,y}^{(i)} t\) comes for all \(t \in F_i\) and \(i \leq n\) for
which \(h t^{-1}\) lies in \(F\).  All in all, the coefficient is at least
\begin{equation}
  \label{eq:13}
  k_n (s_n^2 - 1) - \sum_{t \in F_i, i \leq n} l_n n_{t,y}^{(i)} t \geq k_n (s_n^2 -
  1) - \sum_{t \in F_i, i \leq n} l_n s_n t_n t \geq 0
\end{equation}
by \ref{item:9} from Proposition~\ref{prop:3}.

Now we check that each of the above generators of \(\setQ D_{n,F}\) has a positive
coordinate in \(\setQ D_{n+1,F'}\) which coordinate is zero for the other
generators. This coordinate is \(x_h y1\) for \(x_h y\) where \(h \in F\) and
\(y \in \nSeq\); it is \(x_h c_i^{(n+1)}\) for \(x_h c_i^{(n)}\) where \(h \in F\)
and \(i \leq n\); finally, it is \(v_{n+1,F'}\) for \(v_{n,F}\).

So far we have proved that \(\setQ D_{n,F}\) is an order-subgroup of \(\setQ
D_{n+1,F'}\) if \(F'\) contains \(F\) and \(FF_i\) for \(i \leq n\).  It follows
by induction on \(m-n\) that for every \(n\) and \(F\) and \(m>n\) there is a
finite subset \(S\) of \(H\) such that \(\setQ D_{n,F}\) is an order-subgroup of
\(\setQ D_{m,F'}\) if \(F'\) contains \(S\).  Hence, if \(\setQ D_{n,F}\) is a
subgroup of \(\setQ D_{k,C}\) then both are order-subgroups of \(\setQ D_{m,F'}\) for
suitable \(m\) and \(F'\), hence \(\setQ D_{n,F}\) must be an order-subgroup of
\(\setQ D_{k,C}\).  This shows that the inclusions between the \(\setQ D_{n,F}\) are
order-embeddings.
\end{proof}

Now we return to our second lemma, namely that a cofinal subset of the
\(D_{n,F}\) are order-isomorphic to a finite power of \(\setZ\). 
\begin{proof}[Proof of Lemma~\ref{lem:4}]
If \(n\) is a natural number and \(F\) is a finite subset of \(H\) such that
for all \(h \in H\) the rational numbers \(k_n h^{-1}\) and \(l_n h^{-1}\) are
actually integers then the coefficients of the \(x_h y\) in~\eqref{eq:7} are
integers and hence
\begin{equation}
  \label{eq:11}
  (D_{n,F}, \leq) := (\setZ v_{n,F}, \leq) \times \prod_{\substack{h \in F \\ y \in \nSeq}} (\setZ x_hy, \leq) \times
  \prod_{\substack{i=1 \\ h \in F}}^n (\setZ x_hc_i^{(n)}, \leq).
\end{equation}
We show that such \(D_{n,F}\) form a cofinal system i.e., every \(D_{n,F}\) is
contained in a \(D_{m,F'}\) which has the above property. This is easy once we
know that \(D_{n,F}\) is contained in \(D_{m,F'}\) if \(m \geq n\) and \(F'\)
contains \(F\) and \(FF_i\) for \(i \leq n\). This last statement follows form
the fact that \(x_hc_i^{(n)}\) is contained in \(D_{m,F'}\) if \(m \geq n\) and
\(h\) and \(hF_i\) are contained in \(F'\).  This fact can be proved by
induction on \(m-n\): the case \(m=n\) is obvious because \(h \in F'\).  If \(m >
n\) then \(x_h c_i^{(n+1)}\) is contained in \(D_{m,F'}\) by induction and
\(x_h (c_i^{(n)} - s_n t_n c_i^{(n+1)})\) is also contained in \(D_{m,F'}\) by
Equation~\eqref{eq:12} since \(hF_i\) is contained in \(F'\).
\end{proof}

\section{Automorphisms of dimension groups}
\label{sec:autom-dimens-groups}

In this section we prove that \(H\) is the full group of order-preserving
automorphisms of \(D\) constructed in Proposition~\ref{prop:4}, which finishes
the proof of our \hyperref[th:2]{main theorem}.  This is a special case of the following
proposition, which we are going to prove in this section.
Note that \(D / \setQ u = G\) has the required automorphism group.

\begin{proposition}
  \label{prop:2}
  Let \((D, \leq, u)\) be a dimension group of rank at least \(3\) on which a
  group \(H\) acts by order-preserving automorphisms (the order unit need not
  be preserved).  Let us suppose that the maximal divisible subgroup of \(D\)
  is \(\setQ u\). Furthermore, let us assume that
  \begin{equation}
    \label{eq:4}
    \Aut D / \setQ u = \pm H = \setZ/2\setZ (-1) \times H
  \end{equation}
  i.e., the automorphisms of \(D / \setQ u\) are those induced by \(H\) and their
  negatives. Then \(\Aut (D, \leq) = H\). In other words, all the order-preserving
  automorphisms of \(D\) are those coming from \(H\).
\end{proposition}

We base our proof on the comparison of multiples of \(u\) with elements of
\(D\). This can be described by some rational numbers:
\begin{definition}
  \label{def:1}
  Let \((D, \leq, u)\) be a dimension group. Then for every element \(d\) of
  \(D\) we denote by \(r(d)\) the least rational number \(q\) such that \(q u
  \geq d\). Similarly, let \(l(d)\) denote the greatest rational number \(q\)
  with the property \(q u \leq d\). In other words, for all rational numbers
  \(q\):
  \begin{align}
    \label{eq:20}
    q u \geq d &\iff q \geq r(d),\\
    \label{eq:21}
    q u \leq d &\iff q \leq l(d).
  \end{align}
\end{definition}

We collect the main (and mostly obvious) properties of the functions \(r\) and
\(l\) in the following lemma:
\begin{lemma}
  \label{lem:1}
  Let \((D, \leq, u)\) be a dimension group and \(d\) and element of it. Then the
  following hold:
  \begin{enumerate}[label=(\alph*)]
  \item \label{item:16} The numbers \(r(d)\) and \(l(d)\) exist.
  \item \label{item:17} We have \(l(d)=r(d)=0\) if and only if \(d=0\).
  \item \label{item:19} \(r(-d) = - l(d)\) and \(l(-d)=-r(d)\).
  \item \label{item:22} For all rational numbers \(s\):
    \begin{align}
      r(d+su) &= r(d) + s,\\
      l(d+su) &= l(d) + s.
    \end{align}
  \item \label{item:18} If \(\Phi\) is an order-preserving automorphism of \(D\)
    and \(\Phi(u) = q u\) then
    \begin{align}
      \label{eq:22}
      l(\Phi(d)) &= q l(d),\\
      \label{eq:23}
      r(\Phi(d)) &= q r(d).
    \end{align}
  \item \label{item:20} If \(D\) has rank at least \(3\), the function \(d \mapsto
    l(d) + r(d)\) from \(D\) to the additive group of rational numbers is not
    additive. (It is additive if the rank of \(D\) is at most \(2\).)
  \end{enumerate}
\end{lemma}

\begin{proof}
To prove the existence of \(l(d)\) and \(r(d)\), we may restrict ourselves to
an order-subgroup \((\setZ, \leq)^k\) containing \(d\) and \(u\). Such a subgroup
exists by the definition of dimension group. Clearly, \(u=(n_1,\dotsc,n_k)\)
remains an order unit in the subgroup i.e., its coordinates \(n_i\) are
positive. For every element \((m_1,\dotsc,m_k)\) of \((\setZ, \leq)^k\) the functions
\(r\) and \(l\) are clearly well-defined and have the values
\begin{align}
  \label{eq:24}
  r(m_1,\dotsc,m_k) &= \max_{1 \leq i \leq k} \frac{m_i}{n_i},\\
  \label{eq:25}
  l(m_1,\dotsc,m_k) &= \min_{1 \leq i \leq k} \frac{m_i}{n_i}.
\end{align}
These formulas also show that \(r(d)=l(d)=0\) if and only if \(d=0\). If \(D\)
has rank at least \(3\) then there is an order-subgroup  \((\setZ,
\leq)^k\) of \(D\) containing \(u\) with \(k \geq 3\). We can deduce from the above
formulas that \(r+l\) is not additive even when restricted to such a subgroup.
For example, for the elements \(e_i\) whose \(i\)th coordinate is \(1\) and
all the other coordinates \(0\), we have \(r(e_1)=r(e_2)=r(e_1+e_2)=1\) and
\(l(e_1)=l(e_2)=l(e_1+e_2)=0\), and so \(r(e_1+e_2)+l(e_1+e_2) \neq r(e_1) +
l(e_1) + r(e_2) + l(e_2)\).

The remaining items~\ref{item:19}, \ref{item:22} and \ref{item:18} are obvious.
\end{proof}

Now we start proving the proposition. First we split the order-preserving
automorphism group of \(D\).
\begin{lemma}
  \label{lem:2}
  With the hypothesis and notation of Proposition~\ref{prop:2}, let us denote
  by \(\Gamma\colon \Aut (D, \leq) \to \Aut D / \setQ u\) the canonical map, i.e., \(\Gamma(f)\) is
  the automorphism induced by \(f\) on the quotient. Then there is a
  semidirect product decomposition
  \begin{equation}
    \label{eq:26}
    \Aut (D, \leq) = \Gamma^{-1}(\setZ/2\setZ) \rtimes H, 
  \end{equation}
where \(\setZ/2\setZ\) is generated by \(-1\).
\end{lemma}

\begin{proof}
Note that \(\setQ u\) is invariant under automorphisms since it is the largest
divisible subgroup, so \(\Gamma\) is well-defined. Note that the composition
\begin{equation}
  \label{eq:27}
  H \longrightarrow \Aut (D, \leq) \xrightarrow{\Gamma} \Aut (D / \setQ u) = \setZ/2\setZ \times H \longrightarrow H
\end{equation}
is the identity, which implies the claimed decomposition as a semidirect
product. Here the first arrow is the inclusion of \(H\) given by the
\(H\)-action on \(D\) and the last arrow is projection onto the second
coordinate.
\end{proof}

Now we show that the first term of the semidirect product is trivial, which
finishes the proof.

To this end, we choose an order-preserving automorphism \(\Phi \in \Gamma^{-1}(\setZ/2\setZ)\)
and show that it is the identity. By the choice of \(\Phi\), there is a number
\(\varepsilon=\pm1\) such that the image of \(\Phi-\varepsilon\) is contained in \(\setQ u\). Moreover,
since \(\setQ u\) is invariant, there is a positive rational number \(q\) such that
\(\Phi(u)=qu\).

Our first task is to show that \(q=1\). Therefore we select a nonzero element
\(d\) in the kernel of \(\Phi - \varepsilon\). Since \(\Phi - \varepsilon\) maps to a \(1\)-rank group
\(\setQ u\) and the rank of \(D\) is greater than \(1\), such an element \(d\)
exists.  Now we use Lemma~\ref{lem:1}~\ref{item:18}.  If \(\varepsilon=1\), we obtain
\(r(d)=qr(d)\) and \(l(d)=ql(d)\) and thus \(q=1\) since \(r(d)\) and \(l(d)\)
are not both zero.  If \(\varepsilon=-1\) then \(r(d)=-ql(d)\) and \(l(d)=-qr(d)\).
Again, since at least one of \(r(d)\) and \(l(d)\) is not zero and \(q\) is
positive, \(q\) must be \(1\).

So far we know that \(\Phi(u)=u\). Let \(d\) be an arbitrary element of
\(D\). Then there is a rational number \(s\) depending on \(d\) such that
\(\Phi(d)=\varepsilon d + su\). Our next task is to determine \(s\).

We apply Lemma~\ref{lem:1} again, but this time item \ref{item:22} of it.  If
\(\varepsilon=1\) then \(r(d) = r(d) + s\) and \(l(d) = l(d) + s\).  We conclude that
\(s=0\) for all \(d\).  In other words, \(\Phi\) is the identity.  If \(\varepsilon=-1\)
then we have \(r(d) = s - l(d)\) and \(l(d) = s - r(d)\).  Thus \(s = r(d) +
l(d)\) for all \(d\), which means that \(\Phi(d) = d - (r(d) + l(d)) u\).  So
\(r+l\) is an additive function contradicting Lemma~\ref{lem:1}~\ref{item:20}.

Hence we have proved that \(\Phi=1\) and this finishes the proof.
\providecommand{\bysame}{\leavevmode\hbox to3em{\hrulefill}\thinspace}
\providecommand{\MR}{\relax\ifhmode\unskip\space\fi MR }
\providecommand{\MRhref}[2]{\href{http://www.ams.org/mathscinet-getitem?mr=#1}{#2}
}
\providecommand{\href}[2]{#2}

\end{document}